\documentstyle[12pt]{article}

\newtheorem{t.}{Theorem}[section]
\newtheorem{d.}[t.]{Definition}
\newtheorem{l.}[t.]{Lemma}
\newtheorem{r.}[t.]{Remark}
\newtheorem{p.}[t.]{Proposition}
\newtheorem{c.}[t.]{Corollary}
\newtheorem{e.}[t.]{Example}
\newfont{\kh}{msbm10}
\begin{document}
\title{Induced Norms}
\title{Generalized Induced Norms}
\author{S. Hejazian, M. Mirzavaziri and M. S. Moslehian}
\date{}
\maketitle
\begin{abstract}
Let $\|.\|$ be a norm on the algebra $M_n$ of all $n\times n$ matrices over {\kh C}. An interesting problem in matrix theory is that "are there two norms $\|.\|_1$ and $\|.\|_2$ on {\kh C}$^n$ such that $\|A\|=\max\{\|Ax\|_{2}: \|x\|_{1}=1\}$ for all $A\in M_n$. We will investigate this problem and its various aspects and will discuss under which conditions $\|.\|_1=\|.\|_2$.
\end{abstract}

\section{Preliminaries}

Throughout the paper $M_n$ denotes the complex algebra of all $n\times n$ matrices $A=[a_{ij}]$ with entries in {\kh C} together with the usual matrix operations. Denote by $\{e_1, e_2, \cdots e_n\}$ the standard basis for {\kh C}$^n$, where $e_i$ has $1$ as its $i$th entry and $0$ elsewhere. We denote by $E_{ij}$ the $n\times n$ matrix with $1$ in the $(i,j)$ entry and $0$ elsewhere.

For $1\leq p\leq \infty$ the norm $\ell_p$ on {\kh C}$^n$ is defined as follows:
$$\ell_p(x)=\ell_p(\displaystyle{\sum_{i=1}^n}x_ie_i)=\left \{ \begin{array}{cc}(\displaystyle{\sum_{i=1}^n}|x_i|^p)^{1/p}&1\leq p<\infty\\ \max\{|x_1|, \cdots, |x_n|\}&p=\infty \end{array}\right .$$
A norm $\|.\|$ on {\kh C}$^n$ is said to be unitarily invariant if $\|x\|=\|Ux\|$ for all unitaries $U$ and all $x\in${\kh C}$^n$.

By an algebra norm (or a matrix norm) we mean a norm $\|.\|$ on $M_n$ such that $\|AB\|\leq\|A\|\|B\|$ for all $A, B\in M_n$. An algebra norm $\|.\|$ on $M_n$ is called unitarily invariant if $\|UAV\|=\|A\|$ for all unitaries $U$ and $V$ and all $A\in M_n$. See [2, Chapter IV] for more information.

\begin{e.} {\rm The norm $\|A\|_{\sigma}=\displaystyle{\sum_{i,j=1}^n}|a_{ij}|$ is an algebra norm, but the norm $\|A\|_{m}=\max\{|a_{i,j}|: 1\leq i,j\leq n\}$ is not an algebra norm, since $\|\left [\begin{array}{cc}1&1\\1&1\end{array}\right ]^2\|_{m}>\|\left [\begin{array}{cc}1&1\\1&1\end{array}\right ]\|_{m}^2$.}\end{e.}

\begin{r.}{\rm It is easy to show that for each norm $\| .\|$ on $M_n$, the scaled norm $\max\{\frac{\| AB\|}{\|A \| \|B\|}: A,B\neq 0\}\|.\|$ is an algebra norm; cf. [1, p.114]}\end{r.}

Let $\|.\|_1$ and $\|.\|_2$ be two norms on {\kh C}$^n$. Then for each $A:(${\kh C}$^n,\|.\|_1)\to (${\kh C}$^n,\|.\|_2)$ we can define $\|A\|=\max\{\|Ax\|_{2}: \|x\|_{1}=1\}$. If $\|.\|_1=\|.\|_2$, then $\|I\|=1$ and there are many examples of $\|.\|_1$ and $\|.\|_2$ such that $\|I\|\neq 1$. This shows that given $\|.\|$ on $M_n$, we cannot deduce in general that there is a norm $\|.\|_1$ on {\kh C}$^n$ with $\|A\|=\max\{\|Ax\|_{1}: \|x\|_{1}=1\}$. Let us recall the concept of g-ind norm as follows:

\begin{d.} {\rm Let $\|.\|_1$ and $\|.\|_2$ be two norms on {\kh C}$^n$. Then the norm $\|.\|_{1,2}$ on $M_n$ defined by $\|A\|_{1,2}=\max\{\|Ax\|_2: \|x\|_1=1\}$ is called the generalized induced (or g-ind) norm via $\| .\|_{1}$ and $\| . \|_{2}$. If $\|.\|_1=\|.\|_2$, then $\|.\|_{1,1}$ is called induced norm.}\end{d.}

\begin{e.}{\rm $\|A\|_C=\max\{\displaystyle{\sum_{i=1}^n}|a_{i,j}|: 1\leq j\leq n\}, \|A\|_R=\max\{\displaystyle{\sum_{j=1}^n}|a_{i,j}|: 1\leq i\leq n\}$ and the spectral norm $\|A\|_S= \max\{\sqrt{\lambda}: \lambda {\rm ~is~ an ~eigenvalue~ of} A^*A\}$ are induced by $\ell_1, \ell_\infty$ and $\ell_2$ (or the Eucleadian norm), respectively.

It is known that the algebra norm $\|A\|=\max\{\|A\|_C, \|A\|_R\}$ is not induced [  ] and it is not hard to show that it is not g-ind too; cf. [1, Corollary 3.2.6]}\end{e.}

We need the following proposition which is a special case of a finite dimensional version of the Hahn-Banach theorem [5, p. 104]:

\begin{p.} Let $\|.\|$ be a norm on {\kh C}$^n$ and $y\in${\kh C}$^n$ be a given vector. There exists a vector $y_\circ\in${\kh C}$^n$ such that $y_{\circ}^*y=\|y\|$ and for all $x\in${\kh C}$^n$, $|y_{\circ}^*x|\leq\|x\|$. {\rm (Throughout, $*$ denotes the transpose) [3, Corollary 5.5.15])}\end{p.}

In this paper we examine the following nice problems:\\
(i) Given a norm $\|.\|$ on $M_n$ is there any class ${\cal A}$ of $M_n$ such that the restriction of the norm $\|.\|$ on ${\cal A}$ is g-ind?\\
(ii) When a g-ind norm is unitarily invariant?\\
(iii) If a given norm $\|.\|$ is g-ind via $\|.\|_1$ and $\|.\|_2$, then is it possible to find $\|.\|_1$ and $\|.\|_2$ explicitly in terms of $\|.\|$?\\
(iv) When two g-ind norms are the same?\\
(v) Is there any characterization of the g-ind norms which are algebra norms?

\section{Main Results}

We begin with some observations on generalized induced norms.

Let $\|.\|_{1,2}$ be a generalized induced norm on $M_n$ obtained via $\|.\|_1$ and $\|.\|_2$. Then $\|E_{ij}\|_{1,2}=\max\{\|E_{ij}x\|_2: \|x\|_1=1\}=\max\{\|x_{j}e_i\|_2: \|(x_1,\cdots,x_j,\cdots,x_n)\|_1=1\}=\alpha_j\|e_i\|_2$, where $\alpha_j=\max\{|x_j|: \|(x_1,\cdots,x_j,\cdots,x_n)\|_1=1\}$. In general, for $x\in${\kh C}$^n$ and $1\leq j\leq n$, if $C_{x,j}\in M_n$ is defined by the operator $C_{x,j}(y)=y_{j}x$ then $\|C_{x,j}\|_{1,2}=\alpha_j\|x\|_2$.

Also if for $x\in${\kh C}$^n$ we define $C_{x}\in M_n$ by $C_{x}=\sum_{j=1}^nC_{x,j}$, then clearly $\|C_x\|_{1,2}=\alpha\|x\|_2$, where $\alpha=\max\{|\sum_{j=1}^n y_j|: \|(y_1,\cdots,y_j,\cdots,y_n)\|_1=1\}$.

Now we give a partial solution to Problem (i) and useful direction toward solving Problem (iii):

\begin{p.} Let $\|.\|$ be an algebra norm on $M_n$. Then $\|.\|$ is a g-ind norm on $\{A\in M_n: \|A\|=\|A^{-1}\|=1\}$.\end{p.}

\noindent{\bf Proof.} Put $\|x\|_1=\max\{\|C_{Ax}\|: \|A\|=1\}, \lambda^{-1}=\max\{|\displaystyle{\sum_{i=1}^n}x_i|: \|x\|_1=1\}$ and $\|x\|_2=\lambda\|C_x\|$.

Then we have $\|C_y\|_{1,2}=\max\{\|C_yx\|_2: \|x\|_1=1\}=\max\{|\displaystyle{\sum_{i=1}^n} x_i|\|y\|_2: \|x\|_1=1\}=\|y\|_2\lambda^{-1}=\|C_y\|$.

It follows that for each $y\in${\kh C}$^n$ there is some $x\in${\kh C}$^n$ such that $\|C_yx\|_2=\|C_y\|\|x\|_1=\|C_y\|\max\{\|C_{D x}\|: \|D\|=1\}.$

Now let $A$ be invertible and $\|A^{-1}\|=\|A\|=1$ and $z=A^{-1}C_yx$. Then $\lambda^{-1}\|Bz\|_2=\lambda^{-1}\|BA^{-1}C_yx\|_2=\lambda^{-1}\|Dx\|_2=\|C_{Dx}\|\leq \frac{1}{\|C_y\|}\|C_yx\|_2=\frac{1}{\|C_y\|}\|Az\|_2.$

Now choose $y$ so that $\|C_y\|=1$. Then $\|C_{Bz}\|\leq\|C_{Az}\|$ for all $B\in M_n$. This implies that $\|C_{Az}\|$ is an upper bound for the set $\{\|C_{Bz}\|: \|B\|=1\}$ and indeed $\|C_{Az}\|=\max\{\|C_{Bz}\|: \|B\|=1\}=\|z\|_1$. It follows that $\|A\|=1=\|C_{A(\frac{z}{\|z\|_1})}\|= \max\{\|C_{Au}\|: \|u\|_1=1\}=\max\{\|Au\|_2: \|u\|_1=1\}=\|A\|_{1,2}.\Box$

Let us now answer Question (ii).

\begin{p.} An induced norm $\|.\|_{1,2}$ is unitarily invariant if and only if so are $\|.\|_1$ and $\|.\|_2$.\end{p.}

{\bf Proof.} Let $U, V$ be unital operators and $A$ be an arbitrary operator on {\kh C}$^n$.\\ Suppose that $\|.\|_1$ and $\|.\|_2$ are unitarily invariant. Then $$\|UAV\|_{1,2}=\displaystyle{\max_{x\neq 0}}\frac{\|UAVx\|_2}{\|x\|_1}=\displaystyle{\max_{x\neq 0}}\frac{\|AVx\|_2}{\|x\|_1}=\displaystyle{\max_{y\neq 0}}\frac{\|Ay\|_2}{\|V^{-1}x\|_1}=\displaystyle{\max_{y\neq 0}}\frac{\|Ay\|_2}{\|y\|_1}=\|A_{1,2}.$$
Conversely, if $\|.\|_{1,2}$ is unitarily invariant, then $\|Ux\|_1=\max\{\|AUx\|_2: \|A\|_{1,2}\leq 1\}=\max\{\|Bx\|_2: \|U^{-1}B\|_{1,2}\leq 1\}=\max\{\|Bx\|_2: \|B\|_{1,2}\leq 1\}=\|x\|_1$ and $\|Ux\|_2=\frac{1}{\alpha}\|C_{Ux}\|=\frac{1}{\alpha}\|UC_x\|=\frac{1}{\alpha}\|UC_x\|=\frac{1}{\alpha}\|C_x\|=\|x\|_2.\Box$

Modifying the proof of Theorem 5.6.18 of [3], we obtain a similar useful result for g-ind norms:

\begin{t.} Let $\|.\|_1, \|.\|_2, \|.\|_3$ and $\|.\|_4$ be four given norms on {\kh C}$^n$ and $$R_{i,j}=\max\{\frac{\|x\|_i}{\|x\|_j}: x\neq 0\}, 1\leq i,j\leq 4.$$
Then $$\max\{\frac{\|A\|_{1,2}}{\|A\|_{3,4}}: A\neq 0\}=R_{2,4}R_{3,1}$$

In particual, $\max\{\frac{\|A\|_{1,1}}{\|A\|_{2,2}}: A\neq 0\}=\max\{\frac{\|A\|_{2,2}}{\|A\|_{1,1}}: A\neq 0\}=R_{1,2}R_{2,1}$.\end{t.}

\noindent{\bf Proof.} Let $A$ be a matrix and $x\neq 0$. Then $\frac{\|Ax\|_2}{\|x\|_1}=\frac{\|Ax\|_2}{\|Ax\|_4}.\frac{\|Ax\|_4}{\|x\|_3}.\frac{\|x\|_3}{\|x\|_1}$. Hence $\|A\|_{1,2}\leq R_{2,4}\|A\|_{3,4}R_{3,1}$. Thus $\frac{\|A\|_{1,2}}{\|A\|_{3,4}}\leq R_{2,4}R_{3,1}.$

There are vectors $y, z$ in {\kh C}$^n$ such that $\|y\|_2=\|z\|_2=1, \|y\|_2=R_{2,4}\|y\|_4$ and $\|z\|_3=R_{3,1}\|z\|_1$. By Proposition 1.15, there exists a vector $z_\circ\in$ {\kh C}$^n$ such that $|z_\circ^*x|\leq \|x\|_3$ and $z_\circ^*z=\|z\|_3$.\\ Put $A_\circ=yz_\circ$. Then $\frac{\|A_\circ z\|_2}{\|z\|_1}=\frac{\|yz^*_\circ z\|_2}{\|z\|_1}=\frac{\|y\|_2\|z\|_3}{\|z\|_1}=\|y\|_2R_{3,1}$. Hence $\|A_\circ\|_{1,2}\geq \frac{\|y\|_2}{\|y\|_4}R_{3,1}\|y\|_4=R_{2,4}.R_{3,1}\|y\|_4$. On the other hand, $\frac{\|A_\circ x\|_4}{\|x\|_3}=\frac{\|yz^*_\circ x\|_4}{\|x\|_3}=\frac{\|y\|_4|z^*_\circ x|}{\|x\|_3}\leq\|y\|_4$. Thus $\|A_\circ\|_{3,4}\leq \|y\|_4$. Hence $\frac{\|A_\circ\|_{1,2}}{\|A_\circ\|_{3,4}}\geq \frac{R_{2,4}R_{3,1}\|y\|_4}{\|y\|_4}=R_{2,4}R_{3,1}.\Box$

\begin{c.} (i) $\|.\|_{1,2}\leq\|.\|_{3,2}$ if and only if $\|.\|_1\geq\|.\|_3$,\\
(ii) $\|.\|_{1,2}\leq\|.\|_{1,4}$ if and only if $\|.\|_2\leq\|.\|_4$.\end{c.}

\noindent{\bf Proof.} (i) $\|.\|_{1,2}\leq\|.\|_{3,2}$ if and only if $\max\{\frac{\|A\|_{1,2}}{\|A\|_{3,2}}: A\neq 0\}=R_{2,2}R_{3,1}\leq 1$ and this if and only if $R_{3,1}\leq 1$ or equivalently $\|.\|_3\leq\|.\|_1$. The proof of (ii) is similar.$\Box$

The following corollary completely answers to Question (iv):

\begin{c.} $\|.\|_{1,2}=\|.\|_{3,4}$ if and only if there exists $\gamma>0$ such that $\|.\|_1=\gamma\|.\|_3$ and $\|.\|_2=\gamma\|.\|_4$.\end{c.}

\noindent{\bf Proof.} If $\|A\|_{1,2}=\|A\|_{3,4},$ then $R_{4,2}R_{1,3}=\max\{\frac{\|A\|_{3,4}}{\|A\|_{1,2}}: A\neq 0\}=1=\max\{\frac{\|A\|_{1,2}}{\|A\|_{3,4}}: A\neq 0\}=R_{2,4}R_{3,1}$. Hence $\max\{\frac{\|x\|_2}{\|x\|_4}: x\neq 0\}=R_{2,4}=\frac{1}{R_{3,1}}=\min\{\frac{\|x\|_1}{\|x\|_3}: x\neq 0\}\leq \max\{\frac{\|x\|_1}{\|x\|_3}: x\neq 0\}=R_{1,3}=\frac{1}{R_{4,2}}=\min\{\frac{\|x\|_2}{\|x\|_4}: x\neq 0\}$. Thus there exists a number $\gamma$ such that $\frac{\|x\|_2}{\|x\|_4}=\gamma=\frac{\|x\|_1}{\|x\|_3}.\Box$

\begin{r.}{\rm It is known that each induced norm $\|.\|$ is minimal in the sense that for any matrix norm $\|.\|$, the inequality $\|.\|\leq\|.\|_{1,1}$ implies that $\|.\|=\|.\|_{1,1}$. But this is not true for g-ind norms in general. For instance, put $\|.\|_\alpha=\ell_{\infty}(.), \|.\|_\beta=2\ell_{2}(.)$ and $\|.\|_\gamma=\ell_{2}(.)$. Then $\|.\|_{\gamma,\beta}\leq\|.\|_{\alpha,\beta}$ but $\|.\|_{\gamma,\beta}\neq\|.\|_{\alpha,\beta}$.}
\end{r.}

The following theorem is one of our main theorems and provide a complete solution for Problem (v):

\begin{t.} Let $\|.\|_1$ and $\|.\|_2$ be two norms on {\kh C}$^n$. Then $\|.\|_{1,2}$ is an algebra norm on $M_n$ if and only if $\|.\|_1\leq\|.\|_2$.\end{t.}

\noindent{\bf Proof.} For each $A$ and $B$ in $M_n$ we have
\[\|ABx\|_{2}\leq\|A\|_{1,2}\|Bx\|_{1}\leq\|A\|_{1,2}\|Bx\|_{2}\leq\|A\|_{1,2}\|B\|_{1,2}\|x\|_{1}.\]
Hence $\|AB\|_{1,2}\leq\|A\|_{1,2}\|B\|_{1,2}.$

Conversely, let $\|.\|_{1,2}$ be an algebra norm. Then for each $A,B\in M_{n}$ we have $\Vert AB\|_{2}\leq\| A\|_{1,2}\| B\|_{1,2}\| x\|_{1}$. Let $B$ be an arbitrary member of $M_{n}$. For  $Bx\neq 0$, take $M$ to be the linear span of $\{Bx\}$ and define $f:(M,\|.\|_1)\to${\kh C} by $f(cBx)=\frac{c\| Bx\|_{1}}{\| Bx\|_{2}}$. By the Hahn-Banach Theorem, there is an $F:(${\kh C}$^n,\|.\|_1)\to${\kh C} with $F|_{M}=f$ and $\| F\|=\| f\|=\max\{|f(cBx)|:\| cBx\|_{1}=1\}=\max\{\frac{|c|\| Bx\|_{1}}{\| Bx\|_{2}}:|c|\| Bx\|_{1}=1\}=\frac{1}{\| Bx\|_{2}}$. Define $A:(${\kh C}$^n, \|.\|_1)\to(${\kh C}$^n,\|.\|_2)$ by $Ay=F(y)Bx$. Then $\| A\|_{1,2}=\max\{\| Ay\|_{2}:\| y\|_{1}=1\}=\max\{|F(y)|\| Bx\|_{2}:\| y\|_{1}=1\}=1$, and $\| ABx\|_{2}=|F(Bx)|\| Bx\|_{2}=|f(Bx)|\| Bx\|_{2}=\frac{\| Bx\|_{1}}{\| Bx\|_{2}}\| Bx\|_{2}=\| Bx\|_{1}$. Thus for all $B$,
\[\| Bx\|_{1}=\| ABx\|_{2}\leq \| A\|_{1,2}\| B\|_{1,2}\| x\|_{1}=\| B\|_{1,2}\| x\|_{1},\]
or 
\[\| Bx\|_{1}\leq \| B\|_{1,2}\| x\|_{1}.\]
Now take $N$ to be the linear span of $\{x\}$ and define $g:(N,\|.\|_1)\to${\kh C} by $g(cx)=\frac{c\| x\|_{1}}{\| x\|_{2}}$. By the Hahn-Banach Theorem, there is a $G:(${\kh C}$^n,\|.\|_1)\to${\kh C} with $G|_{N}=g$ and $\| G\|=\| g\|=\max\{|g(cx)|:\| cx\|_{1}\}=\max\{\frac{|c|\| x\|_{1}}{\| x\|_{2}}:|c|\| x\|_{1}=1\}=\frac{1}{\| x\|_{2}}$. Define $B:(${\kh C}$^n,\|.\|_1)\to (${\kh C}$^n,\|.\|_2)$ by $By=G(y)x$. Then $\| B\|_{1,2}=\max\{\| By\|_{2}:\| y\|_{1}=1\}=\max\{|G(y)|\| x\|_{2}:\| y\|_{1}=1\}=\| x\|_{2}\| G\|=1$, and $\| Bx\|_{1}=|G(x)|\| x\|_{1}=|g(x)|\| x\|_{1}=\frac{\| x\|_{1}}{\| x\|_{2}}\| x\|_{1}=\frac{\| x\|_{1}^{2}}{\| x\|_{2}}$.

So
\[\frac{\| x\|_{1}^{2}}{\| x\|_{2}}=\| Bx\|_{1}\leq\| B\|_{1,2}\| x\|_{1}=\| x\|_{1}.\] 
Thus $\|.\|_1\leq\|.\|_2.\Box$

\begin{p.} Suppose that $\|.\|_{1,2}$ is a g-ind norm and $\lambda>0$. Then the scaled norm $\lambda\|.\|_{1,2}$ is a g-ind algebra norm if and only if $\lambda\geq R_{1,2}$.\end{p.}

{\bf Proof.} Evidently, $\lambda\|.\|_{1,2}=\|.\|_{\|.\|_1,\lambda\|.\|_2}$. If $\|.\|_{3,4}=\lambda\|.\|_{1,2}=\|.\|_{\|.\|_1,\lambda\|.\|_2}$ then Corollary 2.5 implies that there exists $\alpha>0$ such that $\|.\|_3=\alpha\|.\|_1$ and $\|.\|_4=\alpha\lambda\|.\|_2$. Now Theorem 2.7 follows that $\lambda\|.\|_{1,2}=\|.\|_{3,4}$ is an algebra norm if and only if $\alpha\|.\|_1\leq\alpha\lambda\|.\|_2$ or equivalently $R_{1,2}\leq \lambda.\Box$

\begin{p.} Let $\|.\|_1$ and $\|.\|_2$ be two norms on {\kh C}$^n$ and $0\neq\alpha,\beta\in${\kh C}. Define $\|.\|_{\alpha}$ and $\|.\|_{\beta}$ on {\kh C}$^n$ by $\|x\|_{\alpha}=\|\alpha x\|_{1}$ and $\|x\|_{\beta}=\|\beta x\|_{2}$, respectively. Then $\|.\|_{\alpha}$  and $\|.\|_{\beta}$ are two norms on {\kh C}$^n$ and $\|.\|_{\alpha,\beta}=|\frac{\beta}{\alpha}|\|.\|_{1,2}$.\end{p.}

\noindent{\bf Proof.} We have $\|A\|_{\alpha,\beta}=\max\{\|Ax\|_{\beta}:\|x\|_{\alpha}=1\}=\max\{\|\beta Ax\|_{2}:\|\alpha x\|_{1}=1\}=|\frac{\beta}{\alpha}|\max\{\|Ay\|_{2}:\|y\|_{1}=1\}=|\frac{\beta}{\alpha}|\|A\|_{1,2}.\Box$

The preceding proposition leads us ti give the following definition:

\begin{d.}{\rm Let $(\|.\|_1,\|.\|_2$) and $(\|.\|_3,\|.\|_4)$ be two pairs of norms on {\kh C}$^n$. We say that $(\|.\|_1,\|.\|_2)$ is generalized induced congruent (gi-congeruent) to $(\|.\|_3,\|.\|_4)$ and we write $(\|.\|_1,\|.\|_2)\equiv_{gi}(\|.\|_3,\|.\|_4)$ if $\|.\|_{1,2}=\gamma\|.\|_{3,4}$ for some $0<\gamma\in${\kh R}.}\end{d.}

Clearly $\equiv_{gi}$ is an equivalence relation. We denote by $[(\|.\|_1,\|.\|_2)]_{gi}$ the equivalence class of $(\|.\|_1,\|.\|_2)$. Proposition 2.9 shows that for each $0<\alpha,\beta\in${\kh R} we have $(\alpha\|.\|_1,\beta\|.\|_2)\equiv_{gi}(\|.\|_1,\|.\|_2)$. Indeed, we have the following result:

\begin{t.} For each pair $(\|.\|_1,\|.\|_2)$ of norms on {\kh C}$^n$ we have $[(\|.\|_1,\|.\|_2)]_{gi}=\{(\alpha\|.\|_1,\beta\|.\|_2):0<\alpha,\beta\in${\kh R}$\}$.\end{t.}

We can extend the above method to find some other norms on $M_n$ which are not necessarily gi-congruent to a given pair $(\|.\|_1,\|.\|_2)$:

\begin{p.} Let $(\|.\|_1,\|.\|_2)$ be a pair of norms on {\kh C}$^n$ and $K,L\in$$M_n$ be two invertible matrices. Define $\|\|_{K}$ and $\|\|_{L}$ and {\kh C}$^n$ by $\|x\|_{K}=\|Kx\|_{1}$ and $\|x\|_{L}=\|Lx\|_{2}$. Then $\|\|_{K}$ and $\|\|_{L}$ are norms on {\kh C}$^n$ and $\|A\|_{K,L}=\|LAK^{-1}\|_{1,2}$.\end{p.}

\noindent{\bf Proof.} Clear and see also Lemma 3.1 of [4].$\Box$

\begin{r.}{\rm Note that the case $K=\alpha I$ and $L=\beta I$ gives Proposition 2.9.}\end{r.}

\end{document}